\documentclass[a4paper,11pt]{article}

\usepackage{amssymb}
\usepackage{amsmath}
\usepackage{tikz}
\usepackage{pgfplots}
\usepackage{amsthm}

\DeclareMathOperator*{\argmax}{arg\,max}
\newcommand{\cupdot}{\mathbin{\mathaccent\cdot\cup}}
\newcommand{\R}{\mathbb{R}}
\renewcommand{\vec}[1]{\boldsymbol{#1}}
\newcommand{\mat}[1]{\boldsymbol{#1}}
\newcommand{\indset}{\mathcal{D}}

\usepackage{algorithm}
\usepackage[noend]{algpseudocode}
\usepackage{float}
\newfloat{algorithm}{t}{lop}

\begin{document}

\author{H.~Harbrecht\thanks{Departement f\"ur Mathematik und Informatik, Unversit\"at Basel, Schweiz}
\and P.~Zaspel\footnotemark[1]}
\title{On the algebraic construction of sparse multilevel approximations of elliptic tensor product problems}
\maketitle

\begin{abstract}
We consider the solution of elliptic problems on the tensor 
product of two physical domains as e.g.~present in the 
approximation of the solution covariance of elliptic partial 
differential equations with random input. Previous 
sparse approximation approaches used a geometrically 
constructed multilevel hierarchy. Instead, 
we construct this hierarchy for a given discretized problem by 
means of the algebraic multigrid method (AMG). Thereby, we 
are able to apply the sparse grid combination technique to 
problems given on complex geometries and for discretizations 
arising from unstructured grids, which was not feasible before. 
Numerical results show that our algebraic construction exhibits
the same convergence behaviour as the geometric construction, while
being applicable even in black-box type PDE solvers.
\end{abstract}

\section{Introduction}
The solution of elliptic problems on tensor products of a 
polygonally bounded domain $\Omega\subset\R^d$ with
e.g.~$d=2,3$ given by
\[\begin{aligned}
(\Delta \otimes \Delta) u &= f & \mbox{on}\ &\Omega\times\Omega\,,\\
u &= 0 & \mbox{on}\ &\partial(\Omega\times\Omega)\,,
\end{aligned}\]
is an important high-dimensional problem. As an example, this problem 
shows up in the estimation of the output covariance of an elliptic partial
differential equation with random input data that is given on a domain 
$\Omega$, see \cite{Harbrecht2010,Harbrecht2013,Schwab2003a,Schwab2003}
for example. The problem becomes high-dimensional since the dimensionality 
of the elliptic problem on $\Omega$ is doubled. In case of real-world problems 
in $d=3$, we end up solving a six-dimensional problem, which might 
become prohibitively expensive.

Recently, there have been developments to overcome this strong limitation. 
These developments are based on the introduction of a geometrically 
constructed multilevel frame, i.e.~a hierarchy of discretizations of the 
elliptic problem on $\Omega$. The multilevel frame gives rise to a sparse 
approximation with respect to the interaction of the involved domains in 
$\Omega\times\Omega$ \cite{Harbrecht2008,Schwab2003a,Schwab2003}. 
It has been shown that the sparse approximation 
allows to solve the tensor product problem at a computational complexity 
that stays essentially (i.e.~up to a poly-logarithmic factor) proportional to 
the number of degrees of freedom to discretize the domain $\Omega$. 
In a more recent work by one of the authors \cite{Harbrecht2013}, 
it has been shown that the sparse approximation can equivalently be 
replaced by the sparse grid combination technique 
\cite{Bungartz2004,Griebel2014,Griebel1992,Hegland2007}. 
This further reduces the computational work and facilitates the implementation.

However, the currently available geometric construction of the multilevel 
hierarchy imposes limitations on the discretization in context of real-world 
problems. First, the coarsest mesh in the hierarchy of discretizations has 
to fully represent the boundary of the geometry $\Omega$. This either 
limits the types of geometry to consider or the computational efficiency
(in case even the coarsest mesh has to be fine at the boundary). Second, 
the use of a fully unstructured mesh becomes barely possible, since we 
are missing a coarsening strategy for such a mesh.   

This work introduces \textit{algebraically} constructed multilevel 
hierarchies \cite{Griebel1994,Griebel.Oswald:2011,Zaspel2016} 
for the solution of elliptic problems on tensor product domains. While 
previous works \cite{Harbrecht2013,Harbrecht2008} first constructed 
the multilevel hierarchy of meshes and then discretized the 
problem by finite elements, the new approach 
first discretizes the problem on $\Omega$ on the finest (potentially 
unstructured) mesh and then constructs coarser versions of the linear 
system resulting from the fine discretization. The coarser problems 
are generated using algebraic coarsening known from the classical 
\textit{Ruge-St\"uben algebraic multigrid} (AMG) \cite{Ruge1986,Stuben2001}.
The algebraic construction of multilevel hieararchies for frames has been previously discussed in context 
of optimal complexity solvers for elliptic problems in \cite{Zaspel2016}. 
However, it has not been applied in the context of sparse approximation 
yet. Note that, by construction, our new approach allows us to overcome 
both the limitations in presence of complex geometries and the requirements 
on the structure of the mesh. Moreover, it perfectly fits into the
context of black-box type PDE solvers.

As it is well-known, a full theory for algebraic multigrid methods, especially
in the multilevel context and on unstructured grids, is still to be developed.
Nevertheless, this technique is extremely popular 
as solver in real-world applications and, usually, empirically shows the 
same performance as geometric multigrid. This work follows the same 
spirit and focuses on the formal construction and the empirical analysis 
of the resulting numerical method. Thereby, we are able to match the convergence 
results available for geometrically constructed sparse approximations, 
while being able to apply this approach to complex geometries and 
unstructured grids in a black-box fashion.

In Section~\ref{sec:univariateConstruction}, the algebraic multilevel 
construction is outlined. This construction is introduced to the tensor product 
problem with sparse approximation and the sparse grid combination technique 
in Section~\ref{sec:bivariateConstruction}. Section~\ref{sec:implementation} 
briefly discusses the implementation. In Section~\ref{sec:numericalResults}, 
we give a series of numerical examples with empirical error analysis. Finally, 
Section~\ref{sec:conclusions} summarizes this work.

\section{Algebraic multilevel constructions}\label{sec:univariateConstruction}
In our algebraic construction, we aim at replacing classical multilevel 
discretizations for elliptic partial differential equations by a purely 
matrix-based construction. That is, we consider an elliptic partial 
differential equation
\begin{equation}\begin{aligned}\label{eq:univariateProblem}
-\Delta u &= f & \mbox{on}\ &\Omega\\
u &= 0 & \mbox{on}\ &\partial\Omega
\end{aligned}\end{equation}
on a polygonally bounded domain $\Omega\subset\R^d$. This problem 
has been discretized by some method on a discretization level $J$, 
leading to a system of linear equations
\begin{equation}
\label{eq:linearSystemFineLevel}
\mat{A}_J \vec{u}_J = \vec{f}_J\,,
\end{equation}
where $\mat{A}_J \in\R^{N_J\times N_J}$ is an \textit{M-matrix}
and $\vec{u}_J,\vec{f}_J\in\R^{N_j}$. In case of the discretization 
by finite elements, $\mat{A}_J$ corresponds to the stiffness matrix 
and $\vec{f}_J$ is the load vector, obtained by, for example, using the 
mass matrix $\mat{M}_J$ and interpolation. Moreover, we identify each 
variable ${u_{J,i}}$ in $\vec{u}_J=(u_{J,1} \ldots u_{J,N_J})^\top$ by its index 
$i$ and introduce the corresponding index set $\indset_J:=\{1,\ldots ,N_J\}$ 
for discretization level $J$.

\subsection{Multilevel hierarchy of discretized problems}
\label{sec:univariateProblemHierarchyOfDiscretizedProblems}
The objective is to construct from \eqref{eq:linearSystemFineLevel} 
a hierarchy of systems of linear equations
\begin{equation}\label{eq:linearSystemHierarchy}
\mat{A}_j \vec{u}_j = \vec{f}_j\,,\quad j=0,\ldots, J\,,
\end{equation}
which are similar to discretizations on different geometric refinement 
levels. Especially, we intend to do this in a purely matrix-based, 
i.e.~algebraic, way by using coarsening and transfer operators from 
algebraic multigrid (AMG) \cite{Stuben2001}. To this end, we first introduce a construction 
method for a hierarchy of variable sets
\begin{equation}\label{eq:variablesHierarchy}
\indset_0 \subset \indset_1 \subset \ldots \subset \indset_J
\end{equation}
of sizes
$$N_0\leq N_2\leq\ldots\leq N_J\,.$$
\begin{algorithm}
  \caption{Standard coarsening algorithm\label{alg:AMGstdCoarsening} \cite{Trottenberg2001}}
  \begin{algorithmic}[1]
    \Statex
    \Require level $j$
    \Function{AMGstandardCoarsening}{}
        \State $F_j:=\emptyset, \indset_{j-1}:=\emptyset, U_j:=\indset_j$
        \For{$i\in U_j$}
        	\State $\lambda_j(i):= \left|{\mathcal{S}_j(i)}^\top\cap U_j\right| +
        	2\,\left|{\mathcal{S}_j(i)}^\top\cap F_j\right|$
        \EndFor
        \While{$\exists i\ \mbox{s.th.}\ \lambda_j(i)\neq 0$}
        	\State find $i_{\max}:= \argmax_i \lambda_j(i)$
        	\State $\indset_{j-1}:=\indset_{j-1}\cup\{i_{\max}\}$
        	\State $U_j:=U_j\setminus\{i_{\max}\}$
        	\For{$j\in ({\mathcal{S}_j(i)}^\top \cap U_j)$}
        		\State $F_j:=F_j\cup\{j\}$
        		\State $U_j:=U_j\setminus\{j\}$
        	\EndFor
	        \For{$i\in U_j$}
    	    	\State $\lambda_j(i):= \left|{\mathcal{S}_j(i)}^\top\cap U_j\right| +
    	    	2\,\left|{\mathcal{S}_j(i)}^\top\cap F_j\right|$
        	\EndFor
        \EndWhile
		\State \Return $\indset_{j-1}, F_j$        
    \EndFunction
  \end{algorithmic}
\end{algorithm} 

In classical Ruge-St\"uben AMG \cite{Ruge1986,Trottenberg2001}, this is achieved by recursively 
splitting the set of variables $\indset_j$ on level $j$ into 
a set of coarse and fine grid variables
$$
\indset_j = \mathcal{D}_{j-1} \cupdot \mathcal{F}_j\,.
$$
Each fine grid variable is supposed to be in the neighborhood of an 
appropriate amount of strongly coupled coarse grid variables, where 
we define the \textit{neighborhood} of a variable $i\in\indset_j$ by
$$
\mathcal{N}_j(i):=\{i^\prime\in\indset_j: i^\prime\neq i,\ a_{j,i i^\prime}\neq 0\}\,,
$$
where $\mat{A}_{j} = (a_{j,i i^\prime})_{i,i^\prime = 1}^{N_j}$. That is, we 
consider neighborhoods between variables by reinterpreting the system 
matrix $\mat{A}_j$ as the adjacency matrix of a graph with edges between 
nodes for each non-zero matrix entry. Moreover, the set of neighboring 
strongly negatively coupled variables of a variable $i$ is
$$
\mathcal{S}_j(i) := \Big\{i^\prime \in\mathcal{N}_j(i) | -a_{j,i i^\prime} \geq \epsilon_{str} \max_k |a_{j,ik}|\Big\}
$$
with a strength measure $0 < \epsilon_{str} < 1$. The \textit{standard coarsening} 
procedure, cf.~Algorithm~\ref{alg:AMGstdCoarsening}~\cite{Trottenberg2001}, builds an appropriate 
splitting $\mathcal{D}_j = \mathcal{D}_{j-1}\cupdot\mathcal{F}_j$ based on these 
considerations. It also involves the sets ${\mathcal{S}_j(i)}^\top$, which are given by
$$
\mathcal{S}_j(i)^\top := \{i^\prime \in\indset_j: i \in \mathcal{S}_j(i^\prime)\}\,.
$$

In order to define the hierarchy of linear systems \eqref{eq:linearSystemHierarchy}, 
we further need a means to transfer information between two consecutive levels 
$j$ and $j+1$. This is done by prolongation operators $\mat{P}_{j}^{j+1}\in
\R^{N_{j+1}\times N_j}$ and restriction operators $\mat{P}_{j+1}^j\in
\R^{N_{j}\times N_{j+1}}$. Prolongation and restriction are done in a purely 
algebraic way based on AMG. In \textit{standard interpolation} \cite{Trottenberg2001}, which 
is one possible type of algebraic prolongation, data given on a fine grid node 
$i\in \mathcal{F}_j$ are interpolated from the set of interpolatory variables
$$
\mathcal{I}_j(i) := \left(\indset_{j-1}\cap \mathcal{S}_j(i)\right)
	\cap \left(\bigcup_{i^\prime \in \mathcal{F}_j
	\cap \mathcal{S}_j(i)} \left(\indset_{j-1}\cap\mathcal{S}_j(i^\prime)\right)\right)\,.
$$
Thus, it is interpolated from strongly negatively coupled coarse grid points 
and all coarse grid points that are strongly negatively coupled to strongly 
negatively coupled fine grid points. The exact choice of prolongation / 
interpolation weights is known from literature \cite{Trottenberg2001}. Restriction is given 
as the transpose of the prolongation, i.e.~$\mat{P}_{j+1}^j = {\mat{P}_{j}^{j+1}}^\top$.

Finally, we recursively define for $j=J-1,\ldots,0$ the matrices and the 
right-hand sides involved in the hierarchy of linear systems
\eqref{eq:linearSystemHierarchy} as
$$
\mat{A}_j := \mat{P}_{j+1}^{j} \mat{A}_{j+1}\mat{P}_{j}^{j+1}, 
\quad \vec{f}_j := \mat{P}_{j+1}^{j}\vec{f}_{j+1}\,.
$$
In order to achieve optimal complexity in AMG, coarser levels are 
constructed such that the \textit{operator complexity}
$$
C_{\mat{A}} := \sum_j \frac{\eta(\mat{A}_j)}{\eta(\mat{A}_J)}\,,
$$
where $\eta(\mat{A}_J)$ is the number of non-zeros in $\mat{A}_J$, 
stays bounded by some constant independent of $J$. If standard 
interpolation and standard coarsening fail in achieving this, stronger 
or more aggressive versions such as \textit{extended / multi-pass interpolation} 
and \textit{aggressive coarsening} on some levels is 
applied to keep this property \cite{Yang2010}. Unfortunately, to the best of the authors' 
knowledge, there is for now no theory on the decay of the number of 
non-zeros in the coarse grid matrices $\mat{A}_j$ constructed by 
classical Ruge-St\"uben AMG on multiple levels and for general M matrices $\mat{A}_J$. The operator 
complexity is therefore always used as empirical measure 
for coarsening quality.

\subsection{Multilevel frames}\label{sec:multilevelFrames}
Let us note here that the above algebraic construction naturally 
leads to \textit{algebraic multilevel frames}, cf.~\cite{Zaspel2016}, for the 
elliptic problem on $\Omega$. That is, we can replace our
original system of linear equations in \eqref{eq:linearSystemFineLevel} 
by the system
\begin{equation}\label{eq:multilevelFramesSystem}
\mat{A}_{\vec{J}} \vec{u}_{\vec{J}} = \vec{f}_{\vec{J}}
\end{equation}
with
$$
\mat{A}_{\vec{J}} := \left(\begin{array}{ccc}
\mat{A}_{11} & \cdots & \mat{A}_{1J}\\
\vdots & \ddots & \vdots\\
\mat{A}_{J1} & \cdots & \mat{A}_{JJ}
\end{array}\right),
\quad \vec{u}_{\vec{J}} := \left(\begin{array}{c}\vec{u}_0\\ \vdots\\ \vec{u}_J\end{array}\right),
\quad \vec{f}_{\vec{J}} := \left(\begin{array}{c}\vec{f}_0\\ \vdots\\ \vec{f}_J\end{array}\right)
$$
and set
$$
\mat{A}_{j_1 j_2} = \mat{P}_{j_1}^{J} \mat{A}_J \mat{P}_{J}^{j_2}\,.
$$
The diagonal matrices $\mat{A}_{jj}$ are the system 
matrices $\mat{A}_j$ from the previous paragraph. Moreover, we 
have extended the prolongation / restriction to arbitrary levels. This 
is possible by concatenating the corresponding operators. Above, we 
further introduce the multi-index $\vec{j}=(j_1,j_2)$ allowing  
the abbreviated notation
$$
\mat{A}_{\vec{J}} = [ \mat{A}_{\vec{j}} ]_{\|\vec{j}\|_{\ell^\infty}\leq J}, 
\quad \vec{u}_{\vec{J}} = [\vec{u}_j]_{|j|\leq J}, 
\quad \vec{f}_{\vec{J}} = [\vec{f}_j]_{|j|\leq J}.
$$
As in multilevel frame discretizations based on geometric refinements / 
coarsening, cf.~\cite{Harbrecht2008}, the above system of linear 
equations now encodes the full information of the hierarchy of systems in 
\eqref{eq:linearSystemHierarchy}. Especially, it is equivalent to the 
the linear system of equations \eqref{eq:linearSystemFineLevel},
if the BPX-preconditioner is applied, cf.\ \cite{BPX,dahmen,Griebel1993a,oswald}. 

The system matrix in \eqref{eq:multilevelFramesSystem} has a 
large kernel, which can be ignored by using appropriate iterative 
linear solvers. Solutions $\vec{u}_{\vec{J}}$ can be projected back 
to single-level solutions $\vec{u}_J$ by applying the operator
$$
\mathcal{P}_{\vec{J}} = \big[ \mat{P}_0^J, \mat{P}_1^J, \dots, \mat{P}_J^J\big]\,.
$$

In \cite{Zaspel2016}, it has been shown by numerical experiments that the 
application of specific iterative solvers to \eqref{eq:multilevelFramesSystem} 
leads to problem-size independent convergence rates also in
case of algebraically constructed multilevel frames.

\section{Sparse algebraic tensor product approach}\label{sec:bivariateConstruction}
Next, we like to consider elliptic problems on tensor products 
$\Omega\times\Omega$ of the polygonally bounded domain 
$\Omega$. That is, we consider problems of the form
\begin{equation}\begin{aligned}\label{eq:tensorProductProblem}
(\Delta \otimes \Delta) u &= f & \mbox{on}\ &\Omega\times\Omega\,,\\
u &= 0 & \mbox{on}\ &\partial(\Omega\times\Omega)\,.
\end{aligned}\end{equation}
As in Section~\ref{sec:univariateConstruction}, we assume to 
have a discretization (e.g.~by finite elements) for the problem 
on a level $J$ resulting in the system of linear
equations
\begin{equation}\label{eq:bivariateLinearSystem}
(\mat{A}_J \otimes \mat{A}_J) \vec{U}_J = \vec{F}_J\,.
\end{equation}
Here, $\mat{A}_J\in\R^{N_J\times N_J}$ is the system matrix 
from \eqref{eq:linearSystemFineLevel}. The operator $\otimes$ 
is the Kronecker product operator for matrices. For matrices 
$\mat{S}\in \R^{n_1\times n_2}$, $\mat{T}\in \R^{m_1\times m_2}$, 
it computes the Kronecker product
$$\mat{S}\otimes \mat{T} := \left(
\begin{array}{ccc}
s_{11}\mat{T} & \hdots & s_{1 n_2}\mat{T}\\
\vdots & \ddots & \vdots\\
s_{n_1 1}\mat{T} & \hdots & s_{n_1 n_2}\mat{T}\\
\end{array}
\right)\,.
$$
Consequently, $(\mat{A}_J\otimes\mat{A}_J)$ becomes a matrix 
of size ${N_J}^2\times {N_J}^2$. Moreover, $\vec{U}_J, \vec{F}_J
\in\R^{N_J\cdot N_J}$ are the solution and the right-hand side, respectively.

By assuming an underlying $d$-dimensional finite element discretization 
with mesh width $h$ and a multigrid-type linear solver, solving the linear
system in~\eqref{eq:bivariateLinearSystem} would require at least $O\left(h^{-2d}\right)$
operations, in contrast to $O(h^{-d})$ for the problem given by
\eqref{eq:linearSystemFineLevel}. This amount of computational work is
prohibitively large, especially for larger $d$. Therefore, we shall find
a way to reduce the amount of work to solve this problem. Before we do that,
we change the problem discretization to a multilevel discretization,
which is the basis for the subsequent sparse approaches.

\subsection{Multilevel frames for tensor product constructions}
\begin{figure}
\begin{tikzpicture}[scale=0.9]
\filldraw[lightgray] (0,0) rectangle (4,4);
\draw[thin,black] (0,0) grid (4,4);
\foreach \x in {0,1,2,3}
{	\foreach \y in {0,1,2,3}
	    \draw (\x+0.5,\y+0.5) node[anchor=center] {\tiny $\widehat{\,\,\mat{A}_{(\x,\y)}}$};}
\end{tikzpicture}\hfill\begin{tikzpicture}[scale=0.9]
\filldraw[lightgray] (0,0) rectangle (1,4);
\filldraw[lightgray] (1,0) rectangle (2,3);
\filldraw[lightgray] (2,0) rectangle (3,2);
\filldraw[lightgray] (3,0) rectangle (4,1);
\draw[thin,black] (0,0) grid (4,4);
\foreach \x in {0,1,2,3}
{	\foreach \y in {0,1,2,3}
	    \draw (\x+0.5,\y+0.5) node[anchor=center] {\tiny $\widehat{\,\,\mat{A}_{(\x,\y)}}$};}
\end{tikzpicture}\hfill\begin{tikzpicture}[scale=0.9]
\filldraw[lightgray] (0,2) rectangle (1,4);
\filldraw[lightgray] (1,1) rectangle (2,3);
\filldraw[lightgray] (2,0) rectangle (3,2);
\filldraw[lightgray] (3,0) rectangle (4,1);
\draw[thin,black] (0,0) grid (4,4);
\foreach \x in {0,1,2,3}
{	\foreach \y in {0,1,2,3}
	    \draw (\x+0.5,\y+0.5) node[anchor=center] {\tiny $\widehat{\,\,\mat{A}_{(\x,\y)}}$};}
\end{tikzpicture}
\caption{\label{fig:subspaces}For discretization level $J=3$, multilevel 
frames on the full tensor product space require a very densely populated 
system matrix $\widehat{\mat{A}_{\vec{J}}}$ \textit{(left)}, while sparse 
approximation leads to the system matrix $\widetilde{\mat{A}_{\vec{J}}}$ 
\textit{(center)} with smaller size due to fewer active (i.e.~gray) matrix 
subblocks. The sparse grid combination technique \textit{(right)} leads 
to the most efficient approximation.}
\end{figure}

As in Section~\ref{sec:univariateProblemHierarchyOfDiscretizedProblems}, 
we can for $\vec{j} = (j, j^\prime)$ introduce multiple levels of systems of
linear equations
\begin{equation}
(\mat{A}_j \otimes \mat{A}_{j^\prime}) \vec{U}_{\vec{j}} = \vec{F}_{\vec{j}}
\end{equation}
with
$$\mat{A}_j= \mat{P}_{J}^{j}\mat{A}_J\mat{P}_{j}^{J}\quad\mbox{and} \quad 
\vec{F}_{\vec{j}} = \left(\mat{P}_{J}^j \otimes \mat{P}_{J}^{j^\prime}\right) \vec{F}_J
$$
by applying coarsening and the transfer operators of algebraic multigrid. 
Thereby, we obtain 
$$
(\mat{A}_j\otimes\mat{A}_{j^\prime})\in\R^{N_j N_{j^\prime}
\times N_j N_{j^\prime}}\quad\text{and}\quad
\vec{U}_{\vec{j}}, \vec{F}_{\vec{j}} \in\R^{N_j N_{j^\prime}}\,.
$$
$\mat{P}_j^{j^\prime}$ is the prolongation / restriction matrix introduced in 
Section~\ref{sec:univariateProblemHierarchyOfDiscretizedProblems}. 

By extending the solution approach in Section~\ref{sec:multilevelFrames} 
to tensor product problems, we finally obtain the multilevel frame linear 
system for the tensor product problem as
\begin{equation}\label{eq:multilevelTPProblem}
\widehat{\mat{A}_{\vec{J}}} \vec{U}_{\vec{J}} =  \vec{F}_{\vec{J}}\,.
\end{equation}
Here, we have
$$
\widehat{\mat{A}_{\vec{J}}} = [ \mat{A}_{j}\otimes\mat{A}_{j^\prime} ]_{j,j^\prime\leq J}
=: [\widehat{\mat{A}_{\vec{j}}}]_{\|\vec{j}\|_{\ell^\infty}\leq J}\,,
$$
and
$$
\vec{U}_{\vec{J}} = [\vec{U}_{\vec{j}}]_{\|\vec{j}\|_{\ell^\infty}\leq J}, 
\quad\vec{F}_{\vec{J}} = [\vec{F}_{\vec{j}}]_{\|\vec{j}\|_{\ell^\infty}\leq J}\,.$$
Note that we indeed construct frames over the tensor product problems 
instead of constructing a tensor product of frame discretizations, 
cf.~\cite{Harbrecht2008}.

In order to characterize the computational complexity for the solution 
of~\eqref{eq:multilevelTPProblem}, we recall that we assume to have 
a constant operator complexity for the sequence of matrices 
$\mat{A}_j$, i.e.~$\sum \eta(\mat{A}_j) \leq c\, \eta(\mat{A}_J)$. 
Moreover, by definition of the Kronecker product, we have the 
number of non-zeros in each block of $\widehat{\mat{A}_{\vec{J}}}$ 
given by
$$
\eta(\mat{A}_j\otimes \mat{A}_{j^\prime}) = \eta(\mat{A}_j)\eta(\mat{A}_{j^\prime})\,.
$$
From that, we can estimate the total number of non-zeros 
in $\widehat{\mat{A}_{\vec{J}}}$ by
$$
\eta(\widehat{\mat{A}_{\vec{J}}}) 
= \sum_{j,j^\prime\leq J} \eta(\mat{A}_j) \eta(\mat{A}_{j^\prime})
\leq c^2 \eta(\mat{A}_J)^2\,.
$$
This means that the computational work to solve 
\eqref{eq:multilevelTPProblem} is asymptotically identical 
to a solve of \eqref{eq:bivariateLinearSystem}.

\subsection{Sparse tensor product construction}
Solving \eqref{eq:bivariateLinearSystem} or 
\eqref{eq:multilevelTPProblem} would be prohibitively expensive, 
cf.~Figure~\ref{fig:subspaces}. As in the geometric multilevel case,
we assume that the solution of the elliptic problem
\eqref{eq:univariateProblem} on $\Omega$ is $H^s$ regular. Therefore,
the solution of the tensor product problem \eqref{eq:tensorProductProblem} 
becomes $H_{\mbox{\tiny mix}}^s$-regular, see \cite{Schwab2003a}.
This allows to follow, for example, the lines of \cite{Harbrecht2008} 
to introduce a sparse, however now algebraically constructed, 
version of the discretized problem. Instead of using all  
sub-problems for multi-indices $\|\vec{j}\|_{\ell^\infty} \leq J$, 
the sparse approximation is reduced to multi-indices 
$\|\vec{j}\|_{\ell^1}\leq J$. Thereby, we obtain a new system
of linear equations
$$
\widetilde{\mat{A}_{\vec{J}}} \widetilde{\vec{U}_{\vec{J}}} 
= \widetilde{\vec{F}_{\vec{J}}}
$$
with
$$
\widetilde{\mat{A}_{\vec{J}}} := [\widehat{\mat{A}_{\vec{j}}}]_{\|\vec{j}\|_{\ell^1}\leq J}, 
\quad \vec{U}_{\vec{J}} = [\vec{U}_{\vec{j}}]_{\|\vec{j}\|_{\ell^1}\leq J}, 
\quad \vec{F}_{\vec{J}} = [\vec{F}_{\vec{j}}]_{\|\vec{j}\|_{\ell^1}\leq J}\,.
$$
Figure~\ref{fig:subspaces} compares both choices in the 
plots on the left-hand side and the center.

As discussed before, there is not much theory on the size of the levels 
in the algebraic multilevel construction. The only available information 
is the assumed bound on the operator complexity. However, this does 
not give enough information to discuss the expected improvement in 
performance due to the sparse construction. Nevertheless, the bound 
implies a similar scaling of the non-zeros with level $j$ as in the 
geometric multilevel construction. Therefore, we here briefly 
discuss the number of non-zeros in $\widetilde{\mat{A}_{\vec{J}}}$ 
for the geometric construction to give a hint towards the possible performance improvement 
by the algebraic sparse construction.

With this in mind, we follow the previous example of (linear) 
finite elements on a mesh with mesh width $h$. The number 
of non-zero entries for matrix $\mat{A}_j$ is proportional to 
the number of elements and therefore
$$
\eta(\mat{A}_j) = O(2^{d\,j})\,.
$$
Moreover, we have $J = O(|\log h|)$. By evaluating 
$$\eta\left(\widetilde{\mat{A}_{\vec{J}}}\right) 
= \sum_{\|\vec{j}\|_{\ell^1}\leq J} \eta\left(\widehat{\mat{A}_{\vec{j}}}\right),
$$
one can easily verify that the number of non-zeros in the 
system matrix in $\widetilde{\mat{A}_{\vec{J}}}$ is asymptotically
$$
\eta\left(\widetilde{\mat{A}_{\vec{J}}}\right) = O\left(|\log h|\, h^{-d}\right)\,.
$$
{
That is, in case a BPX-type preconditioner \cite{BPX,dahmen,Griebel1993a,oswald} and an optimal approach for the
construction of the sub-problem matrices $\widehat{\mat{A}_{\vec{j}}}$ \cite{Balder1996,Bungartz1997,Zeiser2011} is used,
the computational complexity of the problem on the tensor product 
domain $\Omega\times\Omega$ is (up to a logarithmic factor) reduced to 
the computational complexity of the problem on domain $\Omega$.}

\subsection{Sparse grid combination technique}
\label{sec:CT}
It has been shown \cite{Harbrecht2013} that the previous 
sparse approximation is equivalent to the so-called sparse grid 
combination technique. The latter one requires to solve a set of 
decoupled problems
\begin{equation}\label{eq:CTSubproblems}
\widehat{\mat{A}_{\vec{j}}} \vec{U}_{\vec{j}} = \vec{F}_{\vec{j}},
\quad\mbox{where}\ \|\vec{j}\|_{\ell^1}\in\{J,J-1\}\,.
\end{equation}
These are afterwards combined to a solution
\begin{equation}\label{eq:CTCombination}
\widehat{{\vec{U}_{{\vec{J}}}}} 
= \sum_{\|\vec{j}\|_{\ell^1}=J} (\vec{P}_{j}^{J}\otimes\vec{P}_{j^\prime}^J) \vec{U}_{\vec{j}}\,\,\,\, 
-\,\, \sum_{\|\vec{j}\|_{\ell^1}=J-1} (\vec{P}_{j}^{J}\otimes\vec{P}_{j^\prime}^J) \vec{U}_{\vec{j}}\,.
\end{equation}
On the right-hand side of Figure~\ref{fig:subspaces}, the sub-matrices 
$\widehat{\mat{A}_{\vec{j}}}$ used in this approximation have been marked 
gray. As before, one can easily verify that the total number of non-zeros 
of the matrices in \eqref{eq:CTSubproblems} is asymptotically 
$O\left(|\log h|\, h^{-d}\right)$ for the case of linear finite elements 
on a tetrahedral mesh with mesh width $h$ in $d$ dimensions and a 
geometrically constructed multilevel structure. However, Figure~\eqref{fig:subspaces} 
easily clarifies that the pre-asymptotic number of non-zeros in the matrices 
involved in the combination technique is much smaller than the non-zeros 
in the sparse approximation discussed before.

In terms of computational complexity of the combination technique, 
let us remind that the (approximate) solution of each sub-problem in 
\eqref{eq:CTSubproblems} can be realized by an iterative linear solver 
with matrix-vector products. To be more specific, tensor product versions
of standard iterative solvers can be constructed, by reshaping a given
iterate $\vec{U}_{\vec{j}=(j,j^\prime)}\in\R^{N_j\cdot N_j^\prime}$ (and the
appropriate right-hand side) to a 
matrix of size $N_j \times N_{j^\prime}$. Then, the action of one step
of an iterative solver for matrix
$\widehat{\mat{A}_{\vec{j}}}=\mat{A}_j \otimes \mat{A}_{j^\prime}$
is done by first applying the iterative solver step
for $\mat{A}_j$ to all $N_{j^\prime}$ columns of the reshaped matrix and by
second applying the iterative solver step for $\mat{A}_{j^\prime}$ to
all $N_j$ rows of the reshaped matrix.
Since we have all prolongation and restriction operators from AMG 
at our disposal, we can construct, in the above way, a tensor product
version of algebraic multigrid. Given this solver, we obtain roughly 
problem-size independent convergence for each sub-problem in
\eqref{eq:CTSubproblems}, i.e.~we need $O(N_j N_{j^\prime})$ operations
for each sub-problem. In the geometric setting, we would
again have the relation $N_j = O(2^{dj})$ and thereby $O(2^{d(j+j^\prime)})$ operations per sub-problem. Since it holds
$\|\vec{j}\|_{\ell^1}\in\{J,J-1\}$ and the number of sub-problems is $O(J)$,
we would finally end up with a computational complexity of $O(J 2^{dJ})$ or $O(N_J \log N_J)$.

%
%
%

\section{Implementation}\label{sec:implementation}
In our numerical results, we approximate solutions for tensor product finite
element discretizations of elliptic problems based on the combination technique. To this end, 
we assemble system matrices for a given problem, construct the 
multilevel hierarchies, solve the decoupled, anisotropic problems 
in \eqref{eq:CTSubproblems} and combine the solutions following 
the combination rule \eqref{eq:CTCombination}.

\paragraph{Assembly of system matrices.} The discretization by the 
finite element method is done with the \textit{Matlab PDE Toolbox} 
of \textit{Matlab 2017a}. We use linear finite elements and 
construct meshes with maximum element size \texttt{Hmax}$=2^{-J}$. 
Furthermore, we use the option \texttt{Jiggle} to optimize the 
mesh in quality. The stiffness matrix (incorporating boundary 
conditions) is constructed by using the \textit{Matlab} command 
\texttt{assembleFEMatrices} with option \texttt{nullspace}. In 
a similar way, we extract the mass matrix. Afterwards, both 
matrices and the mesh node coordinates are stored to files.

\paragraph{Construction of the multilevel hierarchy.} From 
within \textit{Matlab} we call a newly implemented code that 
uses the parallel linear solver library \textit{hypre} in version 
2.11.1. This library contains the implementation \textit{BoomerAMG} 
of classical Ruge-St\"uben AMG. The code reads the matrix from 
file and creates the AMG multilevel hierarchy by using \textit{hypre}.
In addition to \textit{standard coarsening} with strength measure
$\epsilon_{str}=0.25$ and \textit{standard interpolation}, we use two passes
of \textit{Jacobi} interpolation \cite{Trottenberg2001} with a truncation
of the Jacobi interpolation with a threshold of $0.001$ for the
two-dimensional problems and $0.01$ for the three-dimensional problem.
All other parameters are kept as the defaults of \textit{BoomerAMG}.
After having created the multigrid hierarchy, the program stores 
the prolongation matrices of all created levels to files. These 
are read by \textit{Matlab}.

\paragraph{Solution of the anisotropic tensor product problems.}
{
Based on the prolongation matrices and the system matrix 
$\mat{A}_J$ on the finest levels, the  decoupled problems in 
\eqref{eq:CTSubproblems} can be set up. As discussed before, 
a tensor product version of AMG is used to solve the systems of linear equations.
In our implementation, we construct the sub-problem operators in \eqref{eq:CTSubproblems}
by individually multiplying the transfer operators between two consecutive levels.}
 
Our tensor product AMG is iterated until the convergence criterion
$$
\|\mat{R}_{\vec{j}}^{it}\|_{\ell^2} / \|\vec{F}_{\vec{j}}\|_{\ell^2} \leq \epsilon_{tol}
$$
is fulfilled, where $\mat{R}_{\vec{j}}^{it}$ is the residual of the 
current iterate $\mat{U}_{\vec{j}}^{it}$ in the solver. Since the 
problems in \eqref{eq:CTSubproblems} completely decouple, 
we can easily parallelize their solution process by a \texttt{parfor} 
loop in \textit{Matlab}. In case an individual problem becomes very expensive,
we further implemented a distributed memory parallelization for the
tensor product AMG based on \textit{Matlab}'s \texttt{distributed}
function. Thereby, we overcome the limitation of a non-existing
multi-core parallelization for sparse matrix-vector products in
\textit{Matlab}. 

\paragraph{Combination of the solutions.} In the combination 
phase, we avoid to prolongate the full partial solutions to the 
finest level $J$. Instead, we randomly chose $N_{eval}$ nodes 
on the product of the finest meshes on $\Omega\times\Omega$. 
On these points, we evaluate the combination formula 
\eqref{eq:CTCombination} and compute the empirical 
error measure
$$
e(\mat{U}_{approx}) = \|\mat{U}_{approx} - \mat{U}_{ref}\|_{\ell^2} / \|\mat{U}_{ref}\|_{\ell^2},
$$
where $\mat{U}_{approx}$ is the approximated solution and 
$\mat{U}_{ref}$ is an appropriately evaluated reference solution. 
Note that we do not multiply the tensor product of the 
prolongation with the solution. However, we follow the ideas from
Section~\ref{sec:CT} for the construction of the tensor product AMG
and apply the prolongations direction-wise. The prolongation for
each sub-problem is also parallelized by a \texttt{parfor} loop.

\section{Numerical results}\label{sec:numericalResults}

In our empirical studies, we consider the numerical solution of the problem
\begin{equation}\label{eq:resultsProblem}
\begin{aligned}
(\Delta \otimes \Delta) u &= f & \mbox{on}\ &\Omega\times\Omega\,,\\
u &= 0 & \mbox{on}\ &\partial(\Omega\times\Omega)\,.
\end{aligned}
\end{equation}
by means of the combination technique based on the algebraic multilevel hierarchy.
Different choices will be made for the domain $\Omega$ and the right-hand side $f$.

\subsection{Analytic example on a disk}
\begin{figure}[t]
\scalebox{1.02}{\raisebox{0.6em}{\scalebox{0.4}{\includegraphics{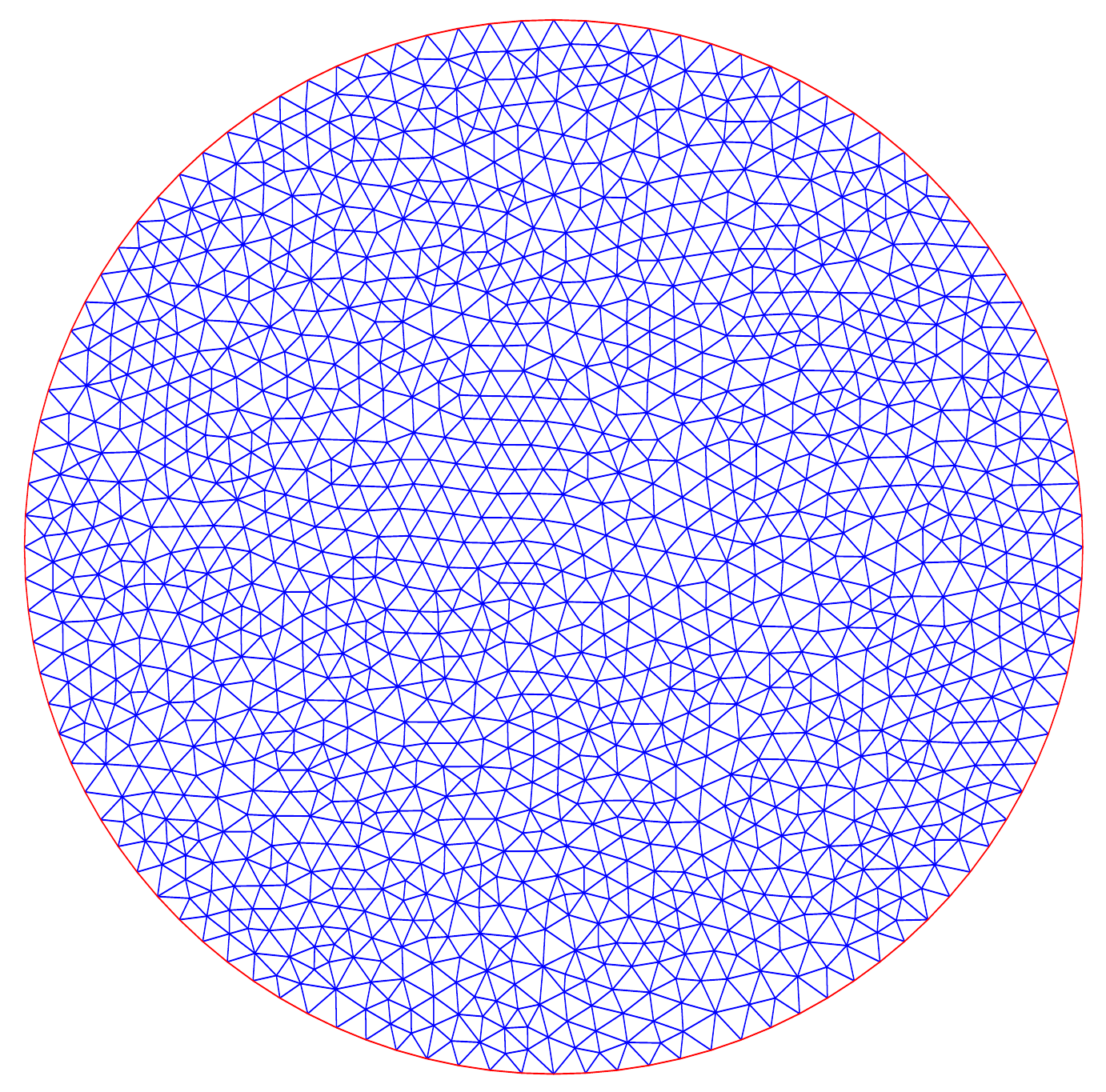}}}\scalebox{0.76}{
\begin{tikzpicture}
        \begin{semilogyaxis}[
          bar width=12pt,
          height=8cm,
          width= 8cm,
          xlabel={level $J$},
          ylabel={relative $\ell^2$ error $e(\vec{U}_{exact})$},
          ylabel near ticks,
          legend style={legend pos = north east, font=\small}
          ]
          \addplot table [x=level, y=error] {study_circle_convergence_new.dat}; \label{plt:circle_error}
          \addlegendentry{error};
	  \addplot[dashed,black] table [x=level, y expr=6*x*4^(-x)]  {study_circle_convergence_new.dat};
	  \addlegendentry{$J 4^{-J}$}
        \end{semilogyaxis}

\end{tikzpicture}}}
\caption{\label{fig:circleResults}The combination technique based on our algebraic multilevel hierarchy and applied to the tensor product of a disk geometry with an unstructured mesh (\textit{left, triangulated with $J=5$}) shows the same convergence as the geometrically constructed combination technique \textit{(right)}.}
\end{figure}
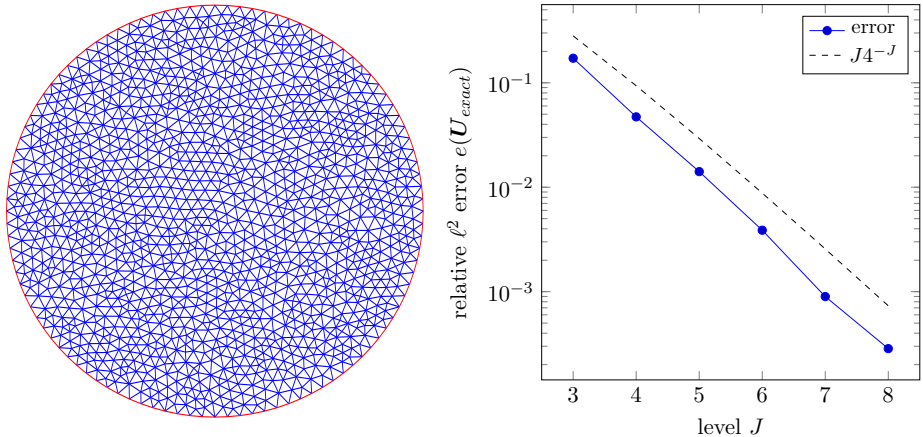

The first study is done on a disk domain $\Omega$ with center $(0,0)^\top$ and radius $0.5$.
We set
$$f(\vec{x},\vec{y}) = 1\,.$$
The exact solution of the resulting problem is
$$u(\vec{x},\vec{y}) = \frac{1}{16} \left(x_1^2 + x_2^2 - 0.5^2\right) \left(y_1^2 + y_2^2 - 0.5^2\right)\,.$$
To approximate the solution $u$ by the combination technique, we follow
the methodology discussed in Section~\ref{sec:implementation}.
As part of this, we triangulate the geometry with a maximum element width of $2^{-J}$.
Figure~\ref{fig:circleResults} shows on the left-hand side the resulting mesh for 
$J=5$. It is obvious that the resulting mesh is unstructured. Therefore, classical
geometric constructions for the sparse grid combination technique would not
be feasible on that mesh. In contrast, our new algebraic approach can solve this problem.

This is shown on the right-hand side of Figure~\ref{fig:circleResults}, where we compare the
numerically approximated solution against the above exact solution.
Convergence results for the choices $J=3,\ldots,8$ are given. From literature,
compare~e.g.~\cite{Harbrecht2013}, we know that the error of the geometrically
constructed sparse grid combination technique scales for the problem under 
consideration like $J 4^{-J}$. As we can see from the convergence results in 
Figure~\ref{fig:circleResults}, the algebraically constructed combination technique 
shows the same convergence behavior, while being applicable to unstructured grids.

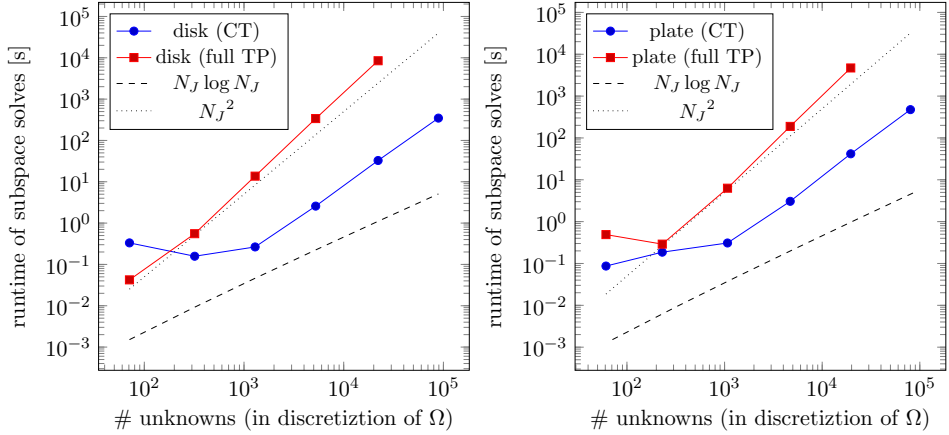
\begin{figure}
\begin{center}
\scalebox{0.76}{
\begin{tikzpicture}
        \begin{loglogaxis}[
          bar width=12pt,
          height=8cm,
          width= 8cm,
          xlabel={\# unknowns (in discretiztion of $\Omega$)},
          ylabel={runtime of subspace solves [s]},
          ylabel near ticks,
          legend style={legend pos = north west, font=\small}
          ]
	  \addplot table [x=N, y=time] {solveTimes_circle_new.dat};
	  \addlegendentry{disk (CT)}
	  \addplot table [x=N, y=time] {solveTimes_circle_full_tp.dat};
	  \addlegendentry{disk (full TP)}
  	  \addplot[dashed] table [x=N, y expr=x*(log2(x)/log2(2.71828))/200000] {solveTimes_circle_new.dat};
	  \addlegendentry{$N_J \log N_J$}
	  \addplot[dotted] table [x=N, y expr=x^2/200000] {solveTimes_circle_new.dat};
	  \addlegendentry{${N_J}^2$}
        \end{loglogaxis}
\end{tikzpicture}}\scalebox{0.76}{
\begin{tikzpicture}
        \begin{loglogaxis}[
          bar width=12pt,
          height=8cm,
          width= 8cm,
          xlabel={\# unknowns (in discretiztion of $\Omega$)},
          ylabel={runtime of subspace solves [s]},
          ylabel near ticks,
          legend style={legend pos = north west, font=\small}
          ]
	  \addplot table [x=N, y=time] {solveTimes_wholeplate_new.dat};
	  \addlegendentry{plate (CT)}
	  \addplot table [x=N, y=time] {solveTimes_wholeplate_full_tp.dat};
	  \addlegendentry{plate (full TP)}
  	  \addplot[dashed] table [x=N, y expr=x*(log2(x)/log2(2.71828))/200000] {solveTimes_circle_new.dat};
	  \addlegendentry{$N_J \log N_J$}
	  \addplot[dotted] table [x=N, y expr=x^2/200000] {solveTimes_wholeplate_new.dat};
	  \addlegendentry{${N_J}^2$}
        \end{loglogaxis}
\end{tikzpicture}}
\end{center}\vspace*{-1em}
\caption{\label{fig:performanceResults}
{We compare the runtime of the new combination technique approach (CT) with the runtime
of the traditional the full tensor-product approach (full TP) for the solution of the tensor product 
elliptic problems on the disk geometry \textit{(left)} and the plate geometry \textit{(right)}}.
}
\end{figure}
\begin{table}
\scalebox{0.73}{
\begin{tabular}{c c|rrrrrrrrrr}
&  & \multicolumn{9}{c}{\bf\# dofs on algebraically coarsened level $\boldsymbol j$}\\
 $\Omega$ & $J$ \textbackslash\ $j$   & 0 & 1  & 2  & 3  & 4  & 5  & 6 &  7 &  8 & 9\\\hline
disk & 3  &  3 & 11 & 28 & 71 &    &    &    &    &   \\
 &  4  & 8  & 21 & 52 & 119 & 320 &    &    &    &   \\
 & 5  & 12  & 31 & 84 & 207 & 495 & 1292&    &    &   \\
 & 6  & 20  & 51 & 139& 348& 852& 2009& 5234&    &   \\
 & 7   & 35 & 93& 244&606& 1473&3510&8415&22118&   \\
   &8   & 46 & 130&366&978&2469&5983&14480&34081&89097\\\hline
plate & 3 & 5  & 20  & 61  &   &   &    &   &    &     \\
      & 4 & 16  & 36  & 90  & 230  &   &    &   &    &     \\
      & 5 & 28  & 68  & 168  & 414  & 1072  &    &   &    &     \\
      & 6 & 46  & 116  &297   & 745  & 1813  & 4703   &   &    &     \\
      & 7 & 63  & 184  & 515  & 1272  & 3117  & 7491   & 19611  &    &     \\
      & 8 & 103  & 302  & 815  & 2124  & 5301  & 12822   & 30639  & 80146   &    \\\hline
spanner & 3 & 4 & 10 & 19 & 50 & 117 & 247 &  &  &  \\  
	& 4 & 11 & 22 & 59 & 147 & 326 & 689 & 1454 & & \\
	& 5 & 40 & 114 & 300 & 689 & 1516 & 3216 & 6484  & 13939 \\
        & 6 & 210 & 548 & 1364 & 3123 & 6708 & 14109 & 29103 & 57438 & 125223\\  
        & 7 & 1386  & 3120 & 6627 & 14016 & 29533 & 61150 & 124921 & 253291 & 496614 & 1082581 \\  
\end{tabular}}
\caption{\label{tab:dofs}For a given problem on level $J$, the algebraic 
multilevel construction on our example domains $\Omega$ constructs 
coarser levels with a decrease of the number of unknowns roughly similar 
to geometric multilevel constructions. Above, only those levels $j$ are 
reported that are used in the convergence study.}
\end{table}

Figure~\ref{fig:performanceResults} shows on the left-hand side computing times
for growing problem size $N_J$ of the univariate discretization of $\Omega$. 
We compare the time required for the solution of the combination technique sub-problems
with the time required to solve the full tensor-product problem \eqref{eq:bivariateLinearSystem} by our tensor-product
AMG implementation. Note that we use the coarse grid hieararchies reported in Table~\ref{tab:dofs} for
both the combination technique and the full tensor-product approach. All measurements were done on a 
compute server with dual 20-core Intel Xeon E5-2698 v4 CPU at 2.2 GHz and 
768 GB RAM. {It becomes evident that our algebraically constructed combination technique
approach beats the full tensor-product approach in both, computational complexity
and effective runtime. However, both results do not show the predicted computational complexity
of $O(N_J \log N_J)$ and $O({N_J}^2)$. There are several reasons for this behavior. 
\begin{itemize}
\item First, algebraic multigrid often shows a small, roughly logarithmic, growth 
in the number of iterations for larger problem sizes, resulting in a slow-down 
by a logarithmic factor. 
\item Second, we observe a certain fill-in 
in the system matrices for coarser problems in the algebraic construction due 
to our choice of an additional (truncated) Jacobi interpolation. However, this 
should be pre-asymptotic behavior. 
\item Third, as can be seen in Table~\ref{tab:dofs}, the AMG 
coarsening approach chosen in our implementation does not show the exact 
same decay rate in the number of levels as we expect from the geometric 
construction. In fact, this leads to a problem-size dependent growth of the coarsest grid.
While this growth does not affect the error decay, it shows up in the computational complexity.
\end{itemize}
Meanwhile,
as stated before, we are able to beat the solution approach based on the full
tensor-product approach in terms of computational complexity. Even more, in terms of runtime,
we are by more than two orders of magnitude faster.}

\subsection{Example on complex geometry with covariance load}
\begin{figure}[t]
\scalebox{1.02}{\raisebox{0.6em}{\scalebox{0.4}{\includegraphics{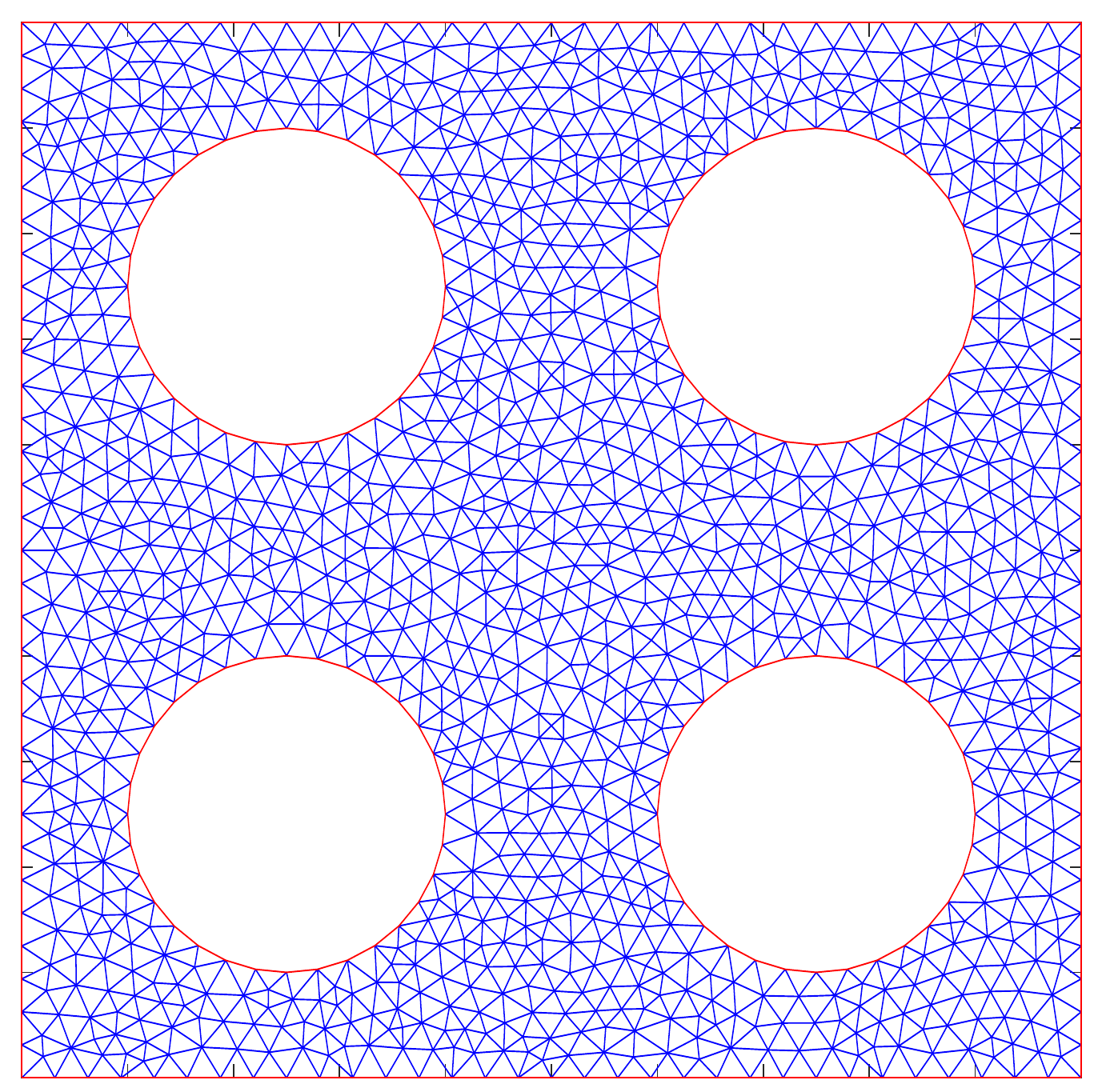}}}\scalebox{0.76}{
\begin{tikzpicture}
        \begin{semilogyaxis}[
          bar width=12pt,
          height=8cm,
          width= 8cm,
          xlabel={level $J$},
          ylabel={relative $\ell^2$ error $e(\vec{U}_{approx})$},
          ylabel near ticks,
          legend style={legend pos = north west, font=\small}
          ]
          \addplot table [x=level, y=error] {study_wholeplate_convergence_new.dat};\label{plt:wholeplate_error}
          \addlegendentry{error};
	  \addplot[dashed,black] table [x=level, y expr=15*x*4^(-x)]  {study_wholeplate_convergence_new.dat};
	  \addlegendentry{$J 4^{-J}$}
        \end{semilogyaxis}

\end{tikzpicture}}}
\caption{\label{fig:holeplateResults}Even for a covariance load on a 
complex geometry (\textit{left}, triangulated for $J=5$), the algebraic 
construction shows the appropriate convergence rate after a short 
pre-asymptotic phase \textit{(right)}.}
\end{figure}
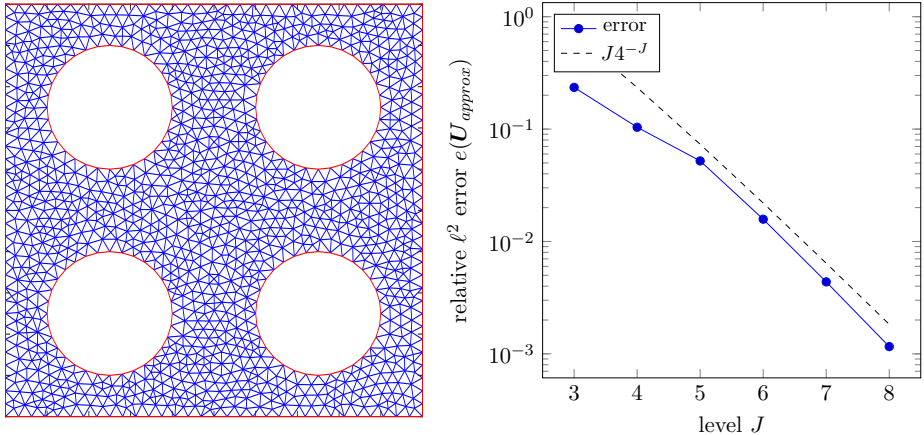

The next numerical study is concerned with the 
solution of the problem \eqref{eq:resultsProblem} with the load
$$f(\vec{x},\vec{y}) = \exp\left(\frac{-\|\vec{x}-\vec{y}\|^2}{\ell}\right)$$
that corresponds to an (unscaled) Gaussian covariance kernel
with correlation length $\ell$. This is a prototype version of the
tensor product elliptic problem on $\Omega\times\Omega$ 
showing up in the computation of the output covariance of 
an elliptic problem on $\Omega$ with random input, 
cf.~\cite{Harbrecht2008}.

In addition to the more complicated right-hand side, we solve the problem
for a rather complex geometry $\Omega$. We choose the geomery of a square
plate on $[0,1]^2$ with circular wholes of radius $0.15$ which are centered at the points
$$\{(0.25, 0.25), (0.25, 0.75), (0.75, 0.25), (0.75,0.75)\}\,.$$
Figure~\ref{fig:holeplateResults} shows its triangulation for $J=5$ on the left-hand side. Note
that it would be almost impossible to solve a problem on such a geometry
with the geometrical construction for the sparse grid combination technique.
However, with the algebraic construction, a coarsening to very few degrees
of freedom becomes easily possible, compare Table~\ref{tab:dofs}.

To be able to compare the above problem against a numerically computed
reference solution, we replace the (sampled) covariance kernel for $\ell=1$
by its low-rank approximation computed with the pivoted Cholesky
factorization \cite{Harbrecht2012}, truncated for a trace norm of $10^{-8}$.
In this case, depending on the problem size,
the truncation results in roughly twenty low-rank terms.

On the right-hand side of Figure~\ref{fig:holeplateResults}, we show the convergence results
with errors computed against the numerically approximated exact solution by
use of the low-rank approximation. After a pre-asymptotic phase, we are able to attain an
error that scales like $J 4^{-J}$ as in the geometric construction.  

The problem size dependent runtime to compute the subspace solutions for the plate geometry
is given in Figure~\ref{fig:performanceResults} {on the right-hand side. We observe similar  
compuational complexities and similar runtimes as in the previous example on the disk.}

\begin{figure}
\begin{center}
\scalebox{0.35}{\includegraphics{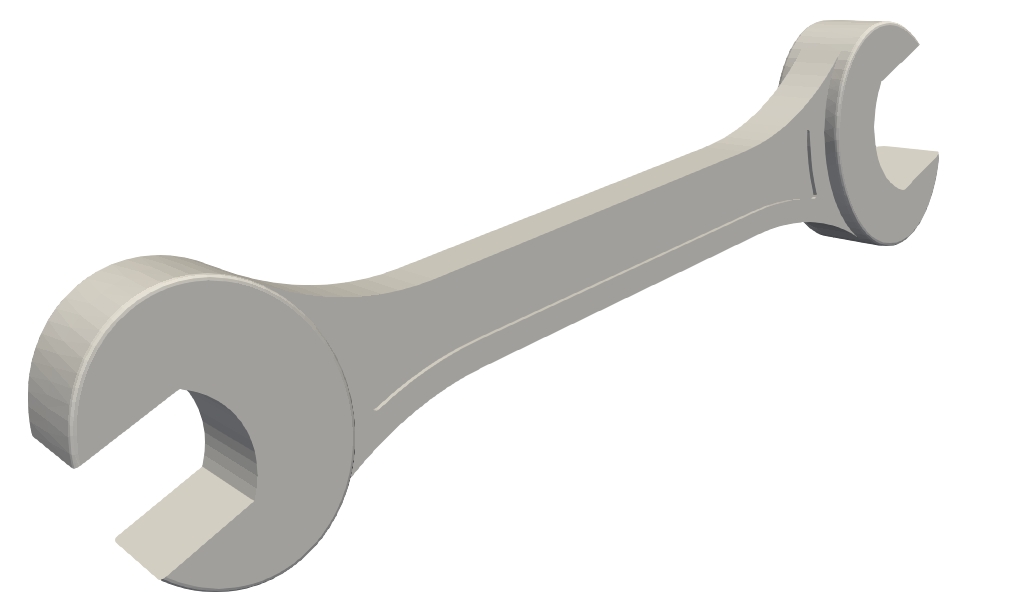}}
\end{center}\vspace*{-1em}
\caption{\label{fig:spannerGeometry}In our large-scale real 
world example, we solve an elliptic problem on the tensor 
product of the three-dimensional geometry of a spanner. For 
a discretization level of $J=7$, the discretization of $\Omega$ 
has more than a million unknowns. This would lead to $10^{12}$, 
that is a \textit{trillion}, unknowns in the full tensor product discretization.}
\end{figure}

\begin{figure}
\begin{center}
\scalebox{0.76}{
\begin{tikzpicture}
        \begin{semilogyaxis}[
          bar width=12pt,
          height=8cm,
          width= 8cm,
          xlabel={level $J$},
          ylabel={relative $\ell^2$ error $e(\vec{U}_{approx})$},
          ylabel near ticks,
          legend style={legend pos = north west, font=\small}
          ]
          \addplot table [x=level, y=error] {study_spanner_convergence.dat};
          \addlegendentry{error};
	  \addplot[dashed,black] table [x=level, y expr=x*6^(-x)*5000]  {study_spanner_convergence.dat};
	  \addlegendentry{$J 6^{-J}$}
        \end{semilogyaxis}

\end{tikzpicture}}
\end{center}\vspace*{-1em}
\caption{\label{fig:spannerResults}{Our algebraic multilevel construction for the sparse grid combination technique on the large-scale three-dimensional spanner geometry gradually approaches the optimal convergence rate of $J 2^{-dJ}$.}}
\end{figure}
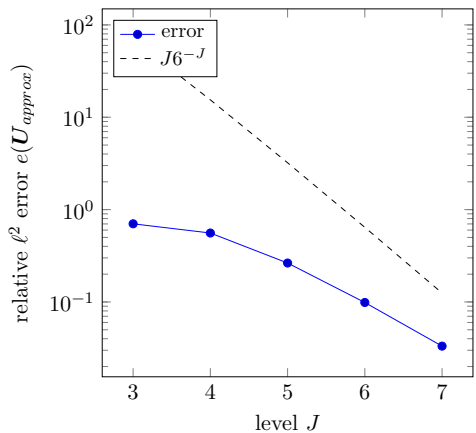

\subsection{Large-scale real-world example}
Our last numerical study treats a large-scale problem with a 
complex real-world geometry $\Omega$. We again aim at solving
\eqref{eq:resultsProblem} for $f(\vec{x},\vec{y}) = 1$. However,
we choose the \textit{three-dimensional} spanner geometry found
in Figure~\ref{fig:spannerGeometry}. In contrast to the previous examples,
we set the maximum mesh width to $2^{5-J}$, since the geometry is
 contained in the rather large bounding box
$[-5,5] \times [-12.2,112] \times [-15.7,15.7]$. Note that the
triangulation of $\Omega$ results for level $J=7$ in a discretization
with $1,082,581$ unknowns. That is, if we would want to solve the full
tensor product problem on $\Omega\times\Omega$, cf.~\eqref{eq:tensorProductProblem},
then we would have to solve a problem with about $10^{12}$,
that is a \textit{trilion}, unknowns. This would be clearly out of scope
even for large parallel clusters. In contrast, the combination
technique allows to solve this problem. Nevertheless, we still have
to solve, e.g.~for level $J=7$ and the system matrix
$\widehat{\,\,\mat{A}_{(0,J)}}$ a problem with $1,082,581 \times 1,386$ unknowns,
compare Table~\ref{tab:dofs}. 

{
In Figure~\ref{fig:spannerResults}, we show the convergence results for this
large scale problem relative to a numerical approximation of the solution. 
Due to the high dimensionality and complexity
of the domain $\Omega$, the convergence results in Figure~\ref{fig:spannerResults} are only gradually
approaching the optimal scaling of $J 2^{-dJ}$. Nevertheless, we are able to solve this problem up to a
certain accuracy. This shows that even very complex problems of
large scale can be solved by the proposed approach.}

\section{Conclusions}\label{sec:conclusions}
In this work, we have introduced an algebraic construction method for 
the sparse approximation of tensor product elliptic problems by means 
of the combination technique. While previous approaches were tight to 
geometric hierarchies of mesh refinements to build the underlying 
multilevel discretization, we were able to solve the given type of 
problems on complex geometries and for unstructured grids by 
an algebraic multilevel hieararchy based on AMG. We could show 
that our approach has the same convergence rates as the geometric 
construction. Measurements of the computational complexity were in 
the linear range with poly-logarithmic factors. Overall, we are now able 
to apply sparse approximation for elliptic tensor product problems in 
a black-box fashion.

\bibliographystyle{abbrv}
\bibliography{bibliography}

\begin{thebibliography}{10}

\bibitem{Balder1996}
R.~Balder and C.~Zenger.
\newblock The solution of multidimensional real {H}elmholtz equations on sparse
  grids.
\newblock {\em SIAM Journal on Scientific Computing}, 17(3):631--646, 1996.

\bibitem{BPX}
J.~Bramble, J.~Pasciak, and J.~Xu.
\newblock {P}arallel multilevel preconditioners.
\newblock {\em Mathematics of Computation}, 55:1--22, 1990.

\bibitem{Bungartz1997}
H.-J. Bungartz.
\newblock A multigrid algorithm for higher order finite elements on sparse
  grids.
\newblock {\em ETNA. Electronic Transactions on Numerical Analysis}, 6:63--77,
  1997.

\bibitem{Bungartz2004}
H.-J. Bungartz and M.~Griebel.
\newblock Sparse grids.
\newblock {\em Acta Numerica}, 13:1--123, 2004.

\bibitem{dahmen}
W.~Dahmen.
\newblock Wavelet and multiscale methods for operator equations.
\newblock {\em Acta Numerica}, 6:55--228, 1997.

\bibitem{Griebel1993a}
M.~Griebel.
\newblock {\em Multilevelmethoden als Iterationsverfahren {\"u}ber
  Erzeugendensystemen}.
\newblock Teubner Skripten zur Numerik. B.G.~Teubner, Stuttgart, 1993.

\bibitem{Griebel1994}
M.~Griebel.
\newblock Multilevel algorithms considered as iterative methods on semidefinite
  systems.
\newblock {\em SIAM International Journal Scientific Statistical Computing},
  15(3):547--565, 1994.

\bibitem{Griebel2014}
M.~Griebel and H.~Harbrecht.
\newblock On the convergence of the combination technique.
\newblock In J.~Garcke and D.~Pfl\"uger, editors, {\em Sparse grids and
  Applications -- Stuttgart 2014}, volume~97 of {\em Lecture Notes in
  Computational Science and Engineering}, pages 55--74. Springer, 2014.

\bibitem{Griebel.Oswald:2011}
M.~Griebel and P.~Oswald.
\newblock Greedy and randomized versions of the multiplicative {S}chwarz
  method.
\newblock {\em Linear Algebra and its Applications}, 7:1596--1610, 2012.

\bibitem{Griebel1992}
M.~Griebel, M.~Schneider, and C.~Zenger.
\newblock A combination technique for the solution of sparse grid problems.
\newblock In P.~de~Groen and R.~Beauwens, editors, {\em {Iterative Methods in
  Linear Algebra}}, pages 263--281. IMACS, Elsevier, North Holland, 1992.

\bibitem{Harbrecht2010}
H.~Harbrecht.
\newblock A finite element method for elliptic problems with stochastic input
  data.
\newblock {\em Applied Numerical Mathematics}, 60(3):227--244, Mar. 2010.

\bibitem{Harbrecht2012}
H.~Harbrecht, M.~Peters, and R.~Schneider.
\newblock On the low-rank approximation by the pivoted {C}holesky
  decomposition.
\newblock {\em Applied Numerical Mathematics}, 62(4):428--440, 2012.

\bibitem{Harbrecht2013}
H.~Harbrecht, M.~Peters, and M.~Siebenmorgen.
\newblock Combination technique based $k$-th moment analysis of elliptic
  problems with random diffusion.
\newblock {\em Journal of Computational Physics}, 252(C):128--141, Nov. 2013.

\bibitem{Harbrecht2008}
H.~Harbrecht, R.~Schneider, and C.~Schwab.
\newblock Multilevel frames for sparse tensor product spaces.
\newblock {\em Numerische Mathematik}, 110(2):199--220, July 2008.

\bibitem{Hegland2007}
M.~Hegland, J.~Garcke, and V.~Challis.
\newblock The combination technique and some generalisations.
\newblock {\em Linear Algebra and its Applications}, 420(2):249--275, 2007.

\bibitem{oswald}
P.~Oswald.
\newblock {\em Multilevel finite element approximation. Theory and
  applications}.
\newblock Teubner Skripten zur Numerik. B.G.~Teubner, Stuttgart, 1994.

\bibitem{Ruge1986}
J.~Ruge and K.~St\"uben.
\newblock Algebraic multigrid ({AMG}).
\newblock In S.~McCormick, editor, {\em Multigrid Methods, Frontiers in Applied
  Mathematics}, volume~5. SIAM, Philadelphia, 1986.

\bibitem{Schwab2003a}
C.~Schwab and R.~A. Todor.
\newblock Sparse finite elements for elliptic problems with stochastic loading.
\newblock {\em Numerische Mathematik}, 95(4):707--734, 2003.

\bibitem{Schwab2003}
C.~Schwab and R.~A. Todor.
\newblock Sparse finite elements for stochastic elliptic problems: Higher order
  moments.
\newblock {\em Computing}, 71(1):43--63, Sept. 2003.

\bibitem{Stuben2001}
K.~St\"uben.
\newblock A review of algebraic multigrid.
\newblock {\em Journal of Computational and Applied Mathematics},
  128(1â2):281--309, 2001.
\newblock {N}umerical Analysis 2000. Vol. VII: Partial Differential Equations.

\bibitem{Trottenberg2001}
U.~Trottenberg and A.~Schuller.
\newblock {\em Multigrid}.
\newblock Academic Press, Inc., Orlando, FL, USA, 2001.

\bibitem{Yang2010}
U.~M. Yang.
\newblock On long-range interpolation operators for aggressive coarsening.
\newblock {\em Numerical Linear Algebra with Applications}, 17(2-3):453--472,
  2010.

\bibitem{Zaspel2016}
P.~Zaspel.
\newblock Subspace correction methods in algebraic multi-level frames.
\newblock {\em Linear Algebra and its Applications}, 488:505--521, 2016.

\bibitem{Zeiser2011}
A.~Zeiser.
\newblock Fast matrix-vector multiplication in the sparse-grid {G}alerkin
  method.
\newblock {\em SIAM Journal of Scientific Computing}, 47(3):328--346, 2011.

\end{thebibliography}

\end{document}